\documentclass[onefignum,onetabnum]{siamart190516}

\usepackage{lipsum}
\usepackage{amsfonts}
\usepackage{graphicx}
\usepackage{epstopdf}
\usepackage{algorithmic}
\usepackage{subcaption}

\usepackage{lmodern}
\usepackage[T1]{fontenc}

\ifpdf
  \DeclareGraphicsExtensions{.eps,.pdf,.png,.jpg}
\else
  \DeclareGraphicsExtensions{.eps}
\fi

\newsiamremark{remark}{Remark}
\newsiamremark{hypothesis}{Hypothesis}
\crefname{hypothesis}{Hypothesis}{Hypotheses}
\newsiamthm{claim}{Claim}

\headers{Faster Stochastic Trace Estimation}{E. Hallman}

\title{Faster Stochastic Trace Estimation with a Chebyshev Product Identity\thanks{
This research was supported in part by the National Science Foundation through grant DMS-1745654.}}

\author{Eric Hallman\thanks{North Carolina State University
  (\email{erhallma@ncsu.edu}, \url{https://erhallma.math.ncsu.edu/}).} }

\usepackage{amsopn}

\definecolor{cerulean}{rgb}{0.16, 0.32, 0.75}

\newcommand{\R}{\mathbb{R}}

\newcommand{\ts}{^T}

\DeclareMathOperator*{\trace}{tr}

\ifpdf
\hypersetup{
  pdftitle={Faster Stochastic Trace Estimation with a Chebyshev Product Identity},
  pdfauthor={Eric Hallman}
}
\fi

\begin{document}

\maketitle

\begin{abstract}
	Methods for stochastic trace estimation often require the repeated evaluation of expressions of the form $z\ts p_n(A)z$, where $A$ is a symmetric matrix and $p_n$ is a degree $n$ polynomial written in the standard or Chebyshev basis. We show how to evaluate these expressions using only $\lceil n/2\rceil$ matrix-vector products, thus substantially reducing the cost of existing trace estimation algorithms that use Chebyshev interpolation or Taylor series. 
	
\end{abstract}

\begin{keywords}
  Chebyshev polynomials, stochastic trace estimation, spectral function, Hutchinson's method
\end{keywords}

\begin{AMS}
  15A63, 68W20, 68W25
\end{AMS}

\section{Introduction}
Given a symmetric matrix $A\in \mathbb{R}^{d\times d}$ and a function $f: \R\rightarrow \R$, we consider the problem of estimating 
\begin{equation}
	\trace f(A) = \sum_{i=1}^d f(\lambda_i), 
\end{equation}
where $\lambda_1,\ldots,\lambda_d$ are the eigenvalues of $A$. When $A$ is large enough to make computing its eigenvalues impractical, one common approach is to use {\it Hutchinson's method} \cite{hutchinson1989stochastic}. This method samples $m$ independent vectors $\{z^{(i)}\}_{i=1}^{m}$ from a Rademacher distribution (entries $\pm 1$ with equal probability), and yields the estimate
\begin{equation}\label{eqn:hutch}
	\trace f(A) \approx \frac{1}{m} \sum_{i=1}^m \left(z^{(i)}\right)\ts f(A)z_i. 
\end{equation}
For many functions $f$ of interest (e.g., $\exp(x)$, $x^{-1}$, $x^{p/2}$, or $\log x$), the right hand side of \eqref{eqn:hutch} is further simplified by approximating $f$ by a degree $n$ polynomial $p_n$, most commonly through Chebyshev interpolation or a Taylor series. 

A common way to evaluate an expression of the form $z\ts p_n(A)z$ is to compute $z_n = p_n(A)z$ and return $z\ts z_n$ (see \cite{han2015large,han2017approximating,boutsidis2017randomized} for examples). We refer to this method as {\it one-sided evaluation}, and it will in general require $n$ matrix-vector products (matvecs) with $A$. For polynomials written in the standard or Chebyshev bases we show how to reduce the number of matvecs to $\lceil n/2\rceil$. Since the cost of many existing trace estimation algorithms is dominated by matvecs, our method will cut their runtime approximately in half. 

\subsection{Related work}
A related idea is explored in \cite{dudley2020monte} where the authors exploit the symmetry of $A$ to reduce the cost of their estimator. Our proposal is more efficient and is more generally applicable, as their method applies only when $f(A)$ is positive definite. 

For more general background on computing matrix polynomials, see \cite[Ch.~4]{higham2008functions} or \cite[Sec.~9.2]{GoVa13}. Our method bears some resemblance to that of Paterson and Stockmeyer \cite{paterson1973number}, but aims to compute $z\ts p_n(A)z$ rather than $p_n(A)$ itself. 

A few recent papers that use Chebyshev approximations for stochastic trace estimation are \cite{han2015large,di2016efficient,han2017approximating,dudley2020monte}, and Taylor series are used similarly in \cite{boutsidis2017randomized}. For applications of stochastic trace estimation, see \cite{ubaru2017applications}. 

One primary competitor to Chebyshev interpolation is stochastic Lanczos quadrature. For more information on this method and its pros and cons with respect to using Chebyshev polynomials, see \cite{ubaru2017fast}. In short, the authors suggest that Lanczos quadrature is generally superior since it converges at twice the rate of Chebyshev interpolation. If this is true, then our methods (which do not extend to Lanczos quadrature) should put Chebyshev interpolation back on more or less equal footing. It may require an approximating polynomial with twice the degree of that needed by Lanczos, but can compute it with the same number of matvecs!

\section{Standard basis}
It was noted in \cite{dudley2020monte} that expressions of the form $z\ts A^nz$ can be evaluated with $\lceil n/2\rceil$ matvecs by letting $k = \lfloor n/2 \rfloor$ and computing $z_k = A^{k}z$, then returning $z_k\ts z_k$ if $n$ is even and $z_k\ts Az_k$ if $n$ is odd. We first extend this idea to polynomials of the form $p_n(x) = \sum_{j=0}^n \alpha_j x^j$. Algorithm \ref{alg:standard} requires $\lceil n/2\rceil$ matvecs and at each step $j$ needs to store only the two most recent vectors $z_j$, $z_{j-1}$ in memory. 

\begin{algorithm}[H]
    \centering
    \caption{Two-sided evaluation (standard basis)}\label{alg:standard}
    \begin{algorithmic}[1]
		\REQUIRE{Symmetric $A\in \R^{d\times d}$, $z_0\in \R^{d}$, polynomial coefficients $a = [\alpha_0,\alpha_1,\ldots,\alpha_n]$}
		\ENSURE{$s = z_0\ts p_n(A)z_0$}
        \STATE{$s = \alpha_0z_0\ts z_0$}
		\FOR{$j = 1,2,\ldots, \lceil n/2\rceil$}
		\STATE{$z_j = Az_{i-1}$}\hfill\COMMENT{$z_j = A^jz_0$}
		\STATE{$s = s + \alpha_{2j-1}z_{j-1}\ts z_i$}
		\STATE{{\bf if} $n = 2j-1$ {\bf then stop}}
		\STATE{$s = s + \alpha_{2j}z_j\ts z_j$}
		\ENDFOR
    \end{algorithmic}
\end{algorithm}

\section{Chebyshev basis}
We use Chebyshev polynomials of the first kind, which can be defined by the recurrence 
\begin{equation}\label{eqn:recurrence}
	T_{j+1}(x) = 2xT_j(x) - T_{j-1}(x), 
\end{equation}
where $T_0(x) = 1$ and $T_1(x) = x$. A function $f: [-1,1]\rightarrow \R$ can then be approximated by a polynomial of the form 
\begin{equation}\label{eqn:chebyBasis}
	f(x) \approx p_n(x) = \sum_{j=0}^n \alpha_j T_j(x). 
\end{equation}
The polynomial $p_n$ interpolates $f$ at a set of Chebyshev nodes $\{x_j\}_{j=0}^n$. Several different choices for the nodes are available \cite{trefethen2008gauss,notaris1997interpolatory}, but as one example Trefethen \cite{trefethen2008gauss} uses the nodes 
\begin{equation}\label{eqn:chebyNodes}
	x_j = \cos \frac{j\pi}{n}, \quad 0 \leq j \leq n
\end{equation}
and shows how to quickly compute the coefficients $\{\alpha_j\}_{j=0}^n$ by using an FFT. Our method applies as long as the polynomial in \eqref{eqn:chebyBasis} is expressed in the Chebyshev basis. 

The key idea is to use the fact that Chebyshev polynomials follow the relation \cite{abramowitz1988handbook}
\begin{equation}
	T_j(x)T_k(x) = \frac{1}{2}\left(T_{j+k}(x) + T_{|k-j|}(x)\right), \quad \forall j,k\geq 0. 
\end{equation}
By letting $k = j$ or $k = j+1$ in the above equation and rearranging, it follows that for all $j\geq 0$,
\begin{equation}
	T_{2j}(x) = 2T_j(x)^2 - T_0(x) = 2T_j(x)^2 - 1
\end{equation}
and 
\begin{equation}
	T_{2j+1}(x) = 2T_j(x)T_{j+1}(x) - T_1(x) = 2T_j(x)T_{j+1}(x) - x.
\end{equation}
We can therefore evaluate terms of the form $z\ts T_{2j}(A)z$ by computing $z_j = T_j(A)z$ and returning $2z_j\ts z_j - z\ts z$. Similarly, we can evaluate terms of the form $z\ts T_{2j+1}(A)z$ by computing $z_{j+1} = T_{j+1}(A)z$ and returning $2z_j\ts z_{j+1} - z\ts Az$. 

Our method is presented in Algorithm \ref{alg:chebyshev}. It requires $\lceil n/2\rceil$ matvecs and at each step $j$ needs to store only the two most recent vectors $z_j$, $z_{j-1}$ in memory. It should therefore take about half the time required by one-sided evaluation. 

\begin{algorithm}[H]
    \centering
    \caption{Two-sided evaluation (Chebyshev basis)}\label{alg:chebyshev}
    \begin{algorithmic}[1]
		\REQUIRE{Symmetric $A\in \R^{d\times d}$, $z_0\in \R^{d}$, Chebyshev coefficients $a = [\alpha_0,\alpha_1,\ldots,\alpha_n]$}
		\ENSURE{$s = z_0\ts p_n(A)z_0$}
		\STATE{$z_1 = Az_0$}
        \STATE{$\zeta_0 = z_0\ts z_0$}
		\STATE{$\zeta_1 = z_0\ts z_1$}
		\STATE{$s = \alpha_0\zeta_0 + \alpha_1\zeta_1 + \alpha_2\left(2(z_1\ts z_1) - \zeta_0\right)$}
		\FOR{$j = 2,3,\ldots, \lceil n/2\rceil$}
		\STATE{$z_j = 2(Az_{i-1})-z_{j-2}$}\hfill\COMMENT{$z_j = T_j(A)z_0$}
		\STATE{$s = s + \alpha_{2j-1}\left(2(z_{j-1}\ts z_j)-\zeta_1\right)$}\label{line:1}
		\STATE{{\bf if} $n = 2j-1$ {\bf then stop}}
		\STATE{$s = s + \alpha_{2j}\left(2(z_j\ts z_j)-\zeta_0\right)$}\label{line:2}
		\ENDFOR
    \end{algorithmic}
\end{algorithm}

\subsection{Stability}
It is possible that due to cancellation and roundoff error, some of the computed values $2z_{j-1}\ts z_j - \zeta_1$ or $2z_{j}\ts z_j - \zeta_0$ may have poor relative accuracy. The {\it absolute} error, however, will be small (on the order of machine precision) compared to $2\|z_{j-1}\|_2\|z_j\|_2 + |\zeta_1|$ and $2\|z_{j}\|_2\|z_j\|_2 + |\zeta_0|$, respectively. Since the first two terms in our sum are $\alpha_0\zeta_0$ and $\alpha_1\zeta_1$, we expect that the effect of these rounding errors will typically be minor. A more rigorous analysis is left for future exploration. 

\section{Numerical Experiments}
We tested our algorithm on synthetic data in order to verify its correctness and to examine how closely the output of Algorithm \ref{alg:chebyshev} aligned with that of one-sided evaluation. We formed a random symmetric matrix $A\in \R^{5000\times 5000}$ and scaled it so that its largest and smallest eigenvalues were $1$ and $-1$. Using the function $f(x) = \exp(10x)$, we computed $\trace f(A)$ to be approximately $7.04\times 10^5$. 

Using the Chebyshev nodes in \eqref{eqn:chebyNodes} with $n=20$, we computed a Chebyshev interpolant $p_n$ to $f$. For each of 100 trials, we drew a random Rademacher vector $z$ and used both our method and one-sided evaluation to compute $z\ts p_n(A)z$. Our method took a total of $5.42$ seconds and one-sided evaluation required $10.70$ seconds, so as predicted two-sided evaluation ran nearly twice as fast. 

Both methods estimated $\trace f(A)$ to be about $7.14\times 10^5$. More notably, the two outputs agreed almost perfectly with each other: the relative error between them was $1.6\times 10^{-16}$, less than machine precision. We then examined the relative errors between the two methods for each evaluation $z\ts p_n(A)z$, and found that the largest was $1.6\times 10^{-14}$. 

To take an even closer look, we examined the values returned by the two algorithms for each individual term $\alpha_jz\ts T_j(A)z$ in each sum. That is, we compared the terms $\alpha_j z_0\ts z_j$ produced by one-sided evaluation to the corresponding terms being added in Lines \ref{line:1} or \ref{line:2} in Algorithm \ref{alg:chebyshev}. The largest relative error found between any single pair of terms was $3.1\times 10^{-11}$, so the two methods of evaluation agreed with each other to a level of accuracy well beyond what would be required for stochastic trace estimation.

\section{Conclusion}
We have shown how to evaluate the expression $z\ts p_n(A)z$, where $A$ is symmetric and $p_n$ is a polynomial in the standard or Chebyshev basis, using no more than $\lceil n/2\rceil$ matvecs with $A$. Our proposed method is simple to implement and can be used for any stochastic trace estimation technique that relies on Taylor expansions or Chebyshev interpolation. The stability of our method remains an open question, but a small numerical experiment suggests that its output will largely agree with that of standard one-sided evaluation. We therefore recommend that two-sided evaluation be incorporated into existing algorithms. 



\bibliographystyle{siamplain}
\bibliography{references}

\begin{thebibliography}{10}

\bibitem{abramowitz1988handbook}
{\sc M.~Abramowitz, I.~A. Stegun, and R.~H. Romer}, {\em Handbook of
  mathematical functions with formulas, graphs, and mathematical tables}, 1988.

\bibitem{boutsidis2017randomized}
{\sc C.~Boutsidis, P.~Drineas, P.~Kambadur, E.-M. Kontopoulou, and A.~Zouzias},
  {\em A randomized algorithm for approximating the log determinant of a
  symmetric positive definite matrix}, Linear Algebra and its Applications, 533
  (2017), pp.~95--117.

\bibitem{di2016efficient}
{\sc E.~Di~Napoli, E.~Polizzi, and Y.~Saad}, {\em Efficient estimation of
  eigenvalue counts in an interval}, Numerical Linear Algebra with
  Applications, 23 (2016), pp.~674--692.

\bibitem{dudley2020monte}
{\sc E.~Dudley, A.~K. Saibaba, and A.~Alexanderian}, {\em Monte carlo
  estimators for the{ S}chatten p-norm of symmetric positive semidefinite
  matrices}, arXiv preprint arXiv:2005.10174,  (2020).

\bibitem{GoVa13}
{\sc G.~H. Golub and C.~F. Van~Loan}, {\em Matrix Computations}, The Johns
  Hopkins University Press, Baltimore, 4th~ed., 2013.

\bibitem{han2017approximating}
{\sc I.~Han, D.~Malioutov, H.~Avron, and J.~Shin}, {\em Approximating spectral
  sums of large-scale matrices using stochastic {C}hebyshev approximations},
  SIAM Journal on Scientific Computing, 39 (2017), pp.~A1558--A1585.

\bibitem{han2015large}
{\sc I.~Han, D.~Malioutov, and J.~Shin}, {\em Large-scale log-determinant
  computation through stochastic {C}hebyshev expansions}, in International
  Conference on Machine Learning, 2015, pp.~908--917.

\bibitem{higham2008functions}
{\sc N.~J. Higham}, {\em Functions of matrices: theory and computation}, SIAM,
  2008.

\bibitem{hutchinson1989stochastic}
{\sc M.~F. Hutchinson}, {\em A stochastic estimator of the trace of the
  influence matrix for {L}aplacian smoothing splines}, Communications in
  Statistics-Simulation and Computation, 18 (1989), pp.~1059--1076.

\bibitem{notaris1997interpolatory}
{\sc S.~E. Notaris}, {\em Interpolatory quadrature formulae with {C}hebyshev
  abscissae of the third or fourth kind}, Journal of computational and applied
  mathematics, 81 (1997), pp.~83--99.

\bibitem{paterson1973number}
{\sc M.~S. Paterson and L.~J. Stockmeyer}, {\em On the number of nonscalar
  multiplications necessary to evaluate polynomials}, SIAM Journal on
  Computing, 2 (1973), pp.~60--66.

\bibitem{trefethen2008gauss}
{\sc L.~N. Trefethen}, {\em Is {G}auss quadrature better than
  {C}lenshaw--{C}urtis?}, SIAM review, 50 (2008), pp.~67--87.

\bibitem{ubaru2017fast}
{\sc S.~Ubaru, J.~Chen, and Y.~Saad}, {\em Fast estimation of tr(f({A})) via
  stochastic {L}anczos quadrature}, SIAM Journal on Matrix Analysis and
  Applications, 38 (2017), pp.~1075--1099.

\bibitem{ubaru2017applications}
{\sc S.~Ubaru and Y.~Saad}, {\em Applications of trace estimation techniques},
  in International Conference on High Performance Computing in Science and
  Engineering, Springer, 2017, pp.~19--33.

\end{thebibliography}
\end{document}